\documentclass[12pt,twoside]{Nreport}
\usepackage{svcon2e}
\usepackage{amsthm,amssymb}
\usepackage{amsbsy,amsfonts,amsmath,amscd}
\usepackage{enumerate}

\catcode`@=12
\makeatletter
\@addtoreset{equation}{section}

\makeatother

\newcommand\Ack{\medskip\noindent{\bf Acknowledgment.}\enspace}

\newcommand{\Proofmark}[1]{} % Otherwise, the section heading gets doubled.

\theoremstyle{plain}
\newtheorem{proposition}{Proposition}[section]
\newtheorem{corollary}{Corollary}[section]
\newtheorem{theorem}{Theorem}[section]
\newtheorem{lemma}{Lemma}[section]

\newtheorem{@definition}{\bf Definition}[section]
\newenvironment{definition}{\begin{@definition}\rm}{\end{@definition}}

\newtheorem{@example}{\bf Example}[section]
\newenvironment{example}{\begin{@example}\rm}{\end{@example}}

\newtheorem{@remark}{\bf Remark}[section]
\newenvironment{remark}{\begin{@remark}\rm}{\end{@remark}}

\newcommand{\F}{\ensuremath{\mathcal{F}}}          % Sheaf F
\renewcommand{\O}{\ensuremath{\mathcal{O}}}        % Structure sheaf

\newcommand{\Z}{\ensuremath{\mathbb{Z}}}           % The integers          
\newcommand{\N}{\ensuremath{\mathbb{N}}}           % The natural numbers
\newcommand{\R}{\ensuremath{\mathbb{R}}}           % The real numbers
\newcommand{\C}{\ensuremath{\mathbb{C}}}           % Complex numbers

\newcommand{\D}{\ensuremath{\mathcal{D}}}          % Differential operators
\newcommand{\End}{\ensuremath{\mathrm{End}}} 
\newcommand{\Aut}{\ensuremath{\mathrm{Aut}}}     
\renewcommand{\L}{\ensuremath{\mathcal{L}}}        % Line bundle
\newcommand{\HC}{\ensuremath{\mathcal{H}}}         % Harish-Chandra modules
\newcommand{\Hom}{\ensuremath{\mathrm{Hom}}}       % Hom
\newcommand{\Cone}{\mathrm{Con}}                   % Mapping cone
\renewcommand{\H}{\ensuremath{\mathrm{H}}}         % Cohomology H
\newcommand{\codim}{\mathrm{codim}}                % codimension

\newcommand{\Dist}{\mathrm{Dist}}                  % Distributions
\newcommand{\Soc}{\mathrm{Soc}}                    % Socle
\newcommand{\hgt}{\mathrm{ht}}                     % Height of a root

\newcommand{\ch}{\mathrm{ch\,}}                    % Lie algebra

\newcommand{\Stab}{\mathrm{Stab}}                  % Lie algebra
\newcommand{\pr}{\mathrm{pr}}                      % Lie algebra
\newcommand{\g}{\mathfrak{g}}                      % Lie algebra
\renewcommand{\b}{\mathfrak{b}}                    % Lie algebra
\newcommand{\p}{\mathfrak{p}}                      % Lie algebra
\newcommand{\h}{\mathfrak{h}}                      % Lie algebra
\newcommand{\n}{\mathfrak{n}}                      % Lie algebra
\newcommand{\spl}{\mathfrak{sl}}
\newcommand{\U}{\mathcal{U}}                       % Env. alg.
\renewcommand{\Im}{\ensuremath{\mathrm{Im}\,}}     % Im(age)
\newcommand{\Coker}{\ensuremath{\mathrm{Coker}}}   % Coker(nel)
\newcommand{\Ext}{\ensuremath{\mathrm{Ext}}}       % Ext
\newcommand{\<}{\langle}                           % Invariant inner product
\renewcommand{\>}{\rangle}

\begin{document}
\pagenumbering{arabic}
\thispagestyle{empty}
\chapter{Twisted Verma modules}
\chapterauthors{H.~H.~Andersen and N.~Lauritzen}

{\renewcommand{\thefootnote}{\fnsymbol{footnote}}

\footnotetext{\kern-19pt{\bf AMS Subject classification:} Primary
; secondary. }

\footnotetext{\kern-19pt{\bf  Keywords and phrases:}
}

\begin{abstract}
Using principal series Harish-Chandra modules, local cohomology 
with support in Schubert cells
and twisting functors we construct certain 
modules parametrized by the Weyl 
group and a highest weight in the subcategory $\O$ of the category of representations
of a complex semisimple Lie algebra. These are in a sense modules between a Verma module and its
dual. We prove that the three different
approaches lead to the same modules. Moreover, we demonstrate that they possess
natural Jantzen type filtrations with corresponding sum formulae.
\end{abstract}

Let $\g$ be a finite dimensional complex semisimple Lie algebra with
a Cartan subalgebra $\h\subset \g$ and Weyl group $W$.
%Fix a
%Borel subalgebra $\b$ and a Cartan subalgebra $\h\subseteq \b$.
In this
paper we consider twisted Verma modules. These are 
in a sense representations between a Verma module and 
its dual. Fix a highest weight $\lambda\in \h^*$. The
twisted Verma modules $M^w(\lambda)$ corresponding to $\lambda$ are
parametrized by the Weyl group $W$. They have the same formal
character as the Verma module $M(\lambda)$ (but in general not the same
module structure). In the affine Kac-Moody setting these 
modules (turning out to be Wakimoto modules) have been studied 
by Feigin and Frenkel \cite{FF}.

We give three rather different ways of constructing twisted
Verma modules. First we obtain them as images of principal
series Harish-Chandra modules (under the Bernstein-Gelfand-Joseph-Enright 
equivalence). In this setting
Irving \cite{Irving} applied wall crossing functors
to describe principal series modules in a regular block inductively 
(Irving uses the term shuffled Verma module for a principal series
module in a regular block). His inductive procedure
inspired this work.

Let $G$ be a complex semisimple algebraic group with Lie algebra $\g$, $B$ a Borel subgroup in $G$, $X=G/B^-$ the flag manifold of $G$, where $B^-$
denotes the Borel subgroup opposite to $B$. 
Let $w_0$ denote the
longest word and $e$ the identity element in the Weyl group $W$ of $G$. We 
let $C(w) = B w B^-/B^- \subseteq X$ denote
the Schubert cell corresponding to $w\in W$. Notice that 
$\codim\, C(w) = \ell(w)$.
It is known that the 
Verma module $M(\lambda)$ with integral highest weight $\lambda$ can
be realized as the top local cohomology group 
$\H^{\ell(w_0)}_{C(w_0)}(X, \L(w_0\cdot\lambda))$ of the line bundle 
$\L(w_0\cdot \lambda)$ with support in the point $C(w_0)$. 
The 
dual Verma module can be realized as the
bottom local cohomology group $\H^0_{C(e)}(X, \L(\lambda))$ of 
the line bundle $\L(\lambda)$ with support in the
big cell $C(e)$.
Our second construction of twisted Verma
modules (which was the starting point of this work) are the intermediate 
local cohomology groups of the line bundle 
$\L(w^{-1}\cdot\lambda)$ with support 
in an arbitrary Bruhat cell $C(w)$  --- these are the modules in 
the global Grothendieck-Cousin complex \cite{Kempf}. The intermediate
local cohomology groups $\H^{\ell(w)}_{C(w)}(X, \L(\lambda))$ are 
isomorphic to dual Verma modules for
dominant weights $\lambda$. In this case the global 
Grothendieck-Cousin complex is
the dual BGG-resolution. 
Let us be more precise about the link from local cohomology to 
principal series modules. Fix a regular antidominant integral weight 
$\lambda$. The principal series modules $M(x, y)$ in the block $\O_\lambda$
(under the Bernstein-Gelfand-Joseph-Enright equivalence) are parametrized
by $(x, y)\in W\times W$. Let 
$C(w) = B w x_0\subseteq X$. Then our result says that
$$
M(x, y) \cong \H^{\ell(x)}_{C(x)}(X, \L(y\cdot\lambda))
$$
as $\g$-modules. Our isomorphism is constructed using
wall translation functors and gives a very explicit algorithm
for obtaining the $\g$-structure of the intermediate local
cohomology groups (starting from a Verma module).
Notice that the local cohomology approach
only makes sense for integral weights.

Following Arkhipov we may for each $w\in W$ define twisting
functors $T_w$ of $\O$ (by tensoring with the ``semiregular''
$U(\g)$-bimodule $S_w$ --- see Section \ref{Sectionsemireg}). When 
applied to a 
Verma module, $T_w$
produces a twisted Verma module. Again it follows quite easily
that the modules obtained in this way satisfy Irving's inductive
procedure. This setup is probably the most powerful for studying
twisted Verma modules and turns out to be the key for showing
that the three approaches are isomorphic: the derived functor
$L T_w$ is
a  self-equivalence of the bounded derived category $D^b(\O)$. 
This implies that twisted Verma modules
only have constant $\g$-endomorphisms (and therefore that they are
indecomposable $\g$-modules). This property allows us to
deduce the required isomorphisms between the three approaches.

The twisting functors also give the required deformation theory for 
constructing Jantzen filtrations and proving sum formulas for
twisted Verma modules (which turn out to be twisted versions of
the original Jantzen sum formula). At the end of the paper we
have used the sum formula to compute the structure
of all twisted Verma modules in the $B_2$-case.

\Ack

We are grateful to S.~Arkhipov for pointing out the 
paper \cite{FF} of Feigin and Frenkel and for 
explaining twisting functors
to us during his stay in Aarhus, January 2001.
We also thank M.~Kashiwara and C.~Stroppel for
discussions that influenced this work. 

\section{Notation and preliminaries}

Fix a {\it complex semisimple Lie algebra\/} $\g$ with a 
{\it Cartan subalgebra\/} $\h$. 
Let $R\subseteq \h^*$ be the {\it root system\/} associated 
with $(\g, \h)$ and $\Z R$ the {\it lattice of roots\/} 
in $\h^*$. Fix a basis $S$ of {\it simple roots\/} and let $R^+$ be 
the {\it positive roots\/} with respect to $S$. Let $\n^+ = \sum_{\alpha>0} 
\g_\alpha$, $\n^- = \sum_{\alpha>0} \g_{-\alpha}$, $\b = \h\oplus \n^+$ and
$\sigma: \g\rightarrow \g$ the {\it Chevalley automorphism\/} of $\g$.
The positive roots give a {\it partial order\/} $\geq$ on $\h^*$ defined by
$\lambda \geq \mu$ if and only if  $\lambda - \mu \in \N R^+$.
The {\it Weyl group\/} $W$ of $R$ acts naturally on $\h^*$.  It is 
generated by 
{\it simple reflections\/} (in the hyperplanes corresponding to the simple roots $S$). Let
$\ell(w)$ denote the {\it length function\/} of an element $w\in W$. We let $e$ and $w_0$
denote the {\it identity element\/} and the {\it unique element of maximal length\/} in $W$ 
respectively. Fix a $W$-invariant positive definite symmetric bilinear
form $(\cdot,\cdot)$ on $\h^*$ and let $\alpha^\vee$ denote the {\it dual
root\/} of $\alpha\in R$ with respect to $(\cdot,\cdot)$.
The {\it dot action\/}
of $W$ on $\h^*$ is given by $w\cdot \lambda = w(\lambda + \rho)-\rho$, where
$2 \rho = \sum_{\alpha>0} \alpha$. A {\it fundamental region\/} for this
action is the set $C$ of {\it antidominant\/} weights translated by $-\rho$. 
Let $\Stab_W(\lambda) \subset W$ denote the {\it stabilizer subgroup\/} of $\lambda\in \h^*$
with respect to the dot-action. A weight $\lambda\in \h^*$ is called
{\it regular\/} if $\Stab_W(\lambda) = \{e\}$ and {\it integral\/} if $\<\lambda, 
\alpha^\vee\>\in \Z$ for every $\alpha\in S$.
We let $C^\circ$ denote the set of regular weights in $C$. 
The {\it enveloping algebra\/} associated with a complex Lie algebra $L$ is denoted
by $U(L)$.

\subsection{The category $\O$}

Let $V$ be 
a (left) $U(\g)$-module.
For $\lambda\in\h^*$, we let $V_\lambda = \{m\in M\mid h m = 
\lambda(h) m, \rm{for\ every\ }h\in\h\}\subseteq V$ denote the {\it weight space\/}
of $V$ corresponding to $\lambda$.   
If $V = \oplus_{\lambda\in \h^*} V_\lambda$, $V$ is called
{\it $\h$-diagonalizable}.

A {\it highest weight vector\/} of weight $\lambda$ in $V$ is a 
non-zero vector in $V_\lambda$ annihilated by $U(\n^+)$. A {\it highest weight
module\/} is a module generated by a highest weight vector.
We let $\O$ (see \cite{BGG}) denote the full subcategory of
the category of left $U(\g)$-modules consisting of modules $V$ 
such that
\begin{itemize}
\item
$V$ is finitely generated
\item
$V$ is $\h$-diagonalizable
\item
$V$ is $U(\n^+)$-finite
\end{itemize}
Any module in $\O$ has a finite filtration with highest weight modules
as subquotients. A highest weight module of weight $\lambda$ is a 
surjective image of the {\it Verma module\/}
$M(\lambda) = U(\g)\otimes_{U(\b)} \C_\lambda\in \O$. Simple $U(\g)$-modules
are parametrized by their highest weight. We let $L(\lambda)$
denote the {\it simple module\/} corresponding to the 
highest weight $\lambda\in \h^*$.
Suppose that $M = \oplus_{\lambda\in \h^*} M_\lambda\in \O$. Then the linear
dual $M^* = \Hom_\C(M, \C)$ is not necessarily $\h$-diagonalizable. This
is remedied by putting $D M = U(\h)$-finite elements in $M^*$. In
fact 
$$
D M = \oplus_{\lambda\in\h^*} M_\lambda^*.
$$
The $\g$-module $D M$ is an object of $\O$ after twisting its
natural $\g$-action with $\sigma$: $X f(m) = f(\sigma(X)m)$, where $X\in \g$ and $f\in D M$. Notice that $D$ then becomes a 
{\it duality\/} 
of $\O$ fixing
simple modules: $D L(\lambda)\cong L(\lambda)$ for $\lambda\in\h^*$. 

We let  
$\ch V = \sum_{\lambda\in \h^*} \dim V_\lambda e^\lambda$ denote
the {\it formal character\/} of an $\h$-diagonalizable $\g$-module $V$ with
finite dimensional weight spaces.

\begin{example}
\label{exampleinO}
Let $V$ be an $\h$-diagonalizabe $\g$-module with finite dimensional weight
spaces. If
$\ch V = \ch M$ for some $M\in \O$, then $V\in \O$.
\end{example}

{\it Composition factors\/} in Verma modules relate to the dot action
by the fundamental result (Harish-Chandra) that $[M(\lambda) : L(\mu)] \neq 0$ implies
that  $\mu \in W\cdot \lambda$.
The category $\O$ decomposes into blocks. We denote for 
$\lambda\in\h^*$ by $\O_\lambda$ the {\it block\/} consisting of 
those $M\in \O$
whose composition factors have the form $L(w\cdot \lambda)$, where $w\in W$. 
Then $\O = \oplus_{\lambda\in C} \O_\lambda$. We let $\pr_\mu$ 
denote the {\it projection\/} $\O \rightarrow \O_\mu$, where $\mu\in \h^*$. 
To a pair of weights $\lambda, \mu\in C$, where $\mu-\lambda$ is integrable 
we have
the
{\it translation functor\/} $T_\lambda^\mu: \O_\lambda\rightarrow \O_\mu$. This functor is
defined by
$$
T_\lambda^\mu(M) = \pr_\mu(M \otimes E), M\in \O_\lambda,
$$
where $E$ is the simple finite dimensional $\g$-module with extremal 
(integral) weight $\mu-\lambda$.
The functors $T_\lambda^\mu$ and $T_\mu^\lambda$ are adjoint.

\begin{definition}
Let $\lambda \in C^\circ$. Pick $\mu\in C$ such that
$\mu\in\lambda + \Z R$ and $\Stab_w(\mu) = \{1, s\}$. This defines the functor
$$
\theta_s=T^\lambda_\mu \circ T_\lambda^\mu : \O_\lambda \rightarrow \O_\lambda
$$
called wall translation (translation through the $s$-wall). Different weights $\mu$ 
with the properties above define naturally isomorphic functors.
The morphism $M\rightarrow \theta_s(M)$ corresponding to
the identity $1\in \Hom(T_\lambda^\mu(M), T_\lambda^\mu(M))$ under
the adjunction isomorphism is called the {\it adjunction morphism\/}.
We let $C_s(M)$ denote the cokernel of the adjunction morphism. 
Thus we have
a short exact sequence
$$
M\rightarrow \theta_s(M) \rightarrow C_s(M)\rightarrow 0.
$$
\end{definition}

\begin{remark} 
The functor $C_s$ is called the shuffling functor in \cite{Irving}.
\end{remark}

\begin{remark}
\label{remarkraphael}
On the level of derived categories the shuffling functor is 
a shadow of the functor 
$$
\tilde{C_s}: X\mapsto \Cone(X \rightarrow \theta_s(X))
$$
where $X\in D^b(\O_\lambda)$, $\Cone$ refers
to the mapping cone of a complex and $\theta_s$ is extended to complexes
in the natural way. If $X$
is concentrated in degree zero and the adjunction morphism
for $X$ is injective, then $\tilde{C_s}(X) \cong C_s(X)$.
Using the fact that $\theta_s^2 \cong \theta_s+\theta_s$ one may
prove that $\tilde{C_s}$ is a self-equivalence of $D^b(\O_\lambda)$. 
We owe
this remark to R.~Rouquier.
\end{remark}

\begin{remark}
For ease of exposition we will restrict ourselves to 
only considering representations with integral weights 
until Section \ref{SectionReform}.
\end{remark}

%\section{The shuffling functor in the derived category}

%In a sense the shuffling functor is not the natural functor to
%study. Instead one should study the functor
%$$
%\tilde{C_s}: X\mapsto \Cone(X \stackrel{\adj}{\rightarrow} \theta_s(X))
%$$
%on the level of complexes $X$ of $\g$-modules in $\O_\lambda$. Suppose
%the complex $X$ is concentrated in degree zero, then 
%$$
%\tilde{C_s}(X) = 0\rightarrow X\stackrel{\adj}{\rightarrow} \theta_s(X) \rightarrow 0.
%$$
%This keeps track not only of the cokernel $C_s(X)$ of the adjunction
%morphism but also the kernel. Also if $\adj(X)$ is injective then
%$\tilde{C_s}(X)$ is resolution of $C_s(X)$, so
%that they agree in the derived category. Let us compute $\tilde{C_s}^2(X)$
%in the derived category, where $X$ is concentrated in degree zero. We need
%to form the mapping cone of
%$$
%\begin{CD}
%
%$$

\section{Formal properties of twisted Verma modules}

In this section we formalize the properties of twisted Verma
modules. We write down a set of properties that characterize
twisted Verma modules in a block $\O_\lambda$ up to
isomorphism.

\subsection{Twisted Verma properties}

\label{twistedvermaproperties}

A family of twisted Verma modules in $\O$ consists of
twisted Verma modules in every block $\O_\lambda$, where
$\lambda\in C$. A set of twisted Verma modules $M_\lambda(x, y)$ in
$\O_\lambda$, where $x, y\in W$ is subject to the following properties
(by abuse of notation we write $M(x, y)$ for
$M_\lambda(x, y)$)

\begin{enumerate}[i.]
\item
\label{itemverma}
$M(e, e)\cong M(\lambda)$.
\item
\label{itemkashiwara}
$M(x, y)\cong M(x s, s y)$ if $x s > x$ and $ s y > y$.
\item 
\label{itemshortexact}
If $\lambda\in C^\circ$
and $y s > y$, where $s$ is a simple reflection, then the adjunction
morphism on $M(x, y)$ is injective and fits the
short exact sequence
$$
0\rightarrow M(x, y)\rightarrow \theta_s M(x, y)\rightarrow 
M(x, ys)\rightarrow 0
$$
for every $x\in W$.
\item
\label{itemwall}
If $\lambda\in C^\circ$ and $\mu\in C$, then 
$$
T_\lambda^\mu M_\lambda(x, y) \cong M_\mu(x, y).
$$
\item
\label{samewall}
If $\lambda\in C^\circ$ then
$$
\theta_s M(x, y) \cong \theta_s M(x, y s)
$$
if $y s > y$.
\end{enumerate}

\begin{lemma}
\label{lemmaverma}
Let $M(x, y)$ be a set of twisted Verma
modules in the block $\O_\lambda$, where $x, y\in W$. 
Then $M(x, e)\cong M(x\cdot\lambda)$
for every $x\in W$ and $\ch M(x, y) = \ch M(x y\cdot\lambda)$ for
all $x, y\in W$.
\end{lemma}
\begin{proof}
By Property \ref{itemwall}) we may reduce to the case, where 
$\lambda\in C^0$.
By Property \ref{itemverma}), $M(e, e)\cong M(\lambda)$. Now
suppose by induction on $\ell(x)$ that $M(x, e)\cong M(x\cdot\lambda)$. Pick
a simple reflection, such that $x s > x$. By Property \ref{itemkashiwara})
it follows that $M(x s, s)\cong M(x, e)$. Also by Property \ref{samewall})
one gets that $\theta_s M(x s, e) \cong\theta_s M(x s, s)$. By Property
\ref{itemshortexact}) $M(x s, e)$ is identified with the kernel of a non-zero
homomorphism $\theta_s M(x\cdot \lambda)\rightarrow M(x\cdot\lambda)$.
This implies that $M(x s, e) \cong M(x s\cdot \lambda)$ 
(\cite{Jantzen1}, 2.17). The fact
that $\ch M(x, y) = \ch M(xy\cdot \lambda)$ follows from 
Property \ref{itemshortexact})
and an easy induction
on $\ell(y)$ ({\it cf.~loc.~cit.\/}).
\end{proof}

\begin{theorem}
\label{theoremunique}
A family of twisted Verma modules is unique up to
isomorphism.
\end{theorem}
\begin{proof}
By Property \ref{itemwall}) it suffices to prove uniqueness
of a set of twisted Verma modules in $\O_\lambda$, where 
$\lambda\in C^\circ$.
Let $y = s_1 \dots s_r$ be a reduced decomposition of $y\in W$, then
$$
M(x, y)\cong C_{s_r}\dots C_{s_1} M(x\cdot 0)
$$
by Lemma \ref{lemmaverma} and Property \ref{itemshortexact}).
\end{proof}

\begin{corollary}
\label{corollaryD}
Let $M(x, y)$ be a set of twisted Verma modules in a block
$\O_\lambda$. Then $D M(x, y) \cong M(x w_0, w_0 y)$.
\end{corollary}
\begin{proof}
We go through the properties for the modules $D M(x w_0, w_0 y)$.
We may assume that $\lambda\in C^0$. Property \ref{itemkashiwara}) implies
that $M(w_0, w_0) \cong M(e, e)$ and hence that 
$M(w_0, w_0)\cong M(\lambda)$ by Property \ref{itemverma}). But $M(\lambda)\cong D M(\lambda)$ as $M(\lambda)$ is simple. To verify
Property \ref{itemkashiwara}), assume that $x s > x$ and $s y > y$. Then
write $s w_0 = w_0 t$ for a suitable simple reflection $t$ and
therefore $D M(x s w_0, w_0 s y) \cong D M(x w_0 t, t w_0 y) \cong
D M( x w_0, w_0 y)$. Suppose that $y s > y$, then Property \ref{itemshortexact})
follows from applying Property \ref{samewall}) and
dualizing the short exact sequence
$$
0\rightarrow M(x w_0, w_0 y s)\rightarrow \theta_s M(x w_0, w_0 y s)
\rightarrow M(x w_0, w_0 y)\rightarrow 0.
$$
Properties \ref{itemwall}) and \ref{samewall}) are immediate
using that translation commutes with duality $D$. By 
Theorem \ref{theoremunique} we get that $D M(x, y)\cong M(x w_0, w_0 y)$.
\end{proof}

\begin{lemma}
\label{lemmaunique}
Suppose that there exists a family of twisted Verma modules
in $\O$ admitting only constant $\g$-endomorphisms. A
family of modules
satisfying all properties
of twisted Verma modules except that we only have
a short exact sequence
$$
0\rightarrow M(x, y)\rightarrow \theta_s M(x, y)\rightarrow 
M(x, ys)\rightarrow 0
$$
in Property \ref{itemshortexact}) (without any conditions on
the morphisms involved) is a family of twisted Verma modules.
\end{lemma}
\begin{proof}
Let $M'(x, y)$ denote the modules in the family of
twisted Verma modules with only constant $\g$-endomorphisms
and $M(x, y)$ the modules in the other family in a regular
block $\O_\lambda$. We will prove that $M'(x, y)\cong M(x, y)$
by induction on $\ell(y)$. As 
$$
\End_\g(M_\mu'(x, y)) = 
\Hom_\g(T_\lambda^\mu M'_\lambda(x, y), T_\lambda^\mu M'_\lambda(x, y)),
$$ 
we get $\Hom_\g(M'_\lambda(x, y), 
\theta_s(M'_\lambda(x, y)) ) = \C$, so that the morphism 
$M(x, y) \rightarrow \theta_s M(x, y)$ in the
relaxed Property \ref{itemshortexact}) has to
be a constant multiple of the adjunction morphism.
By the proof of Lemma \ref{lemmaverma},
$M_1(x, e)\cong M(x, e) \cong M(x\cdot \lambda)$ (the morphism 
$M(x, e)\rightarrow \theta_s M(x, e)$ in the relaxed Property 
\ref{itemshortexact}) 
has to be a constant multiple
of the adjunction morphism). Now suppose that 
$M_1(x, y)\cong M(x, y)$ and let $s$ be a simple reflection with
$y s > y$. Then 
$$
\Hom_\g(M(x, y), \theta_s M(x, y))\cong\Hom_\g(M_1(x, y), \theta_s M_1(x, y))\cong \C
$$ 
and we have a commutative diagram
$$
\begin{CD}
0@>>> M_1(x, y)@>>> \theta_s M_1(x, y)@>>> M_1(x, ys)@>>> 0\\
@.      @|              @|                @.   \\ 
0@>>> M(x, y)  @>>> \theta_s M(x, y)@>>>   M(x, ys)@>>> 0
\end{CD}
$$
giving an isomorphism $M_1(x, ys) \cong M(x, ys)$.
\end{proof}

\begin{remark} \label{remStroppel}
Using Remark \ref{remarkraphael} one may show that a
twisted Verma module in a regular block only has
constant $\g$-endomorphisms. We need this result not only
in regular blocks but in
arbitrary blocks (or at least in (semiregular) blocks $\O_\lambda$, 
where $\lambda$ is
stabilized by at most one simple reflection ). This
is where the twisting functor approach is very useful.
In the semiregular case C.~Stroppel has proved that the principal series 
modules admit only 
constant $\g$-endomorphisms using results of Joseph on 
completion functors. 
\end{remark}

\section{Principal series Harish-Chandra modules}

Here we recall basic properties of and results on
Harish-Chandra modules following \cite{Irving} and
\cite{Jantzen2}. The goal is to prove that principal
series Harish-Chandra modules when viewed in $\O$ through
the categorical equivalence of Bernstein et.~al., 
form a family of twisted Verma modules. Basically
this has been done by Irving \cite{Irving}. Here we
reformulate his results in our setup.

\subsection{Definition}

Let $M$ be a $\g\times\g$-module and view $M$ as a $\g$-module
through the embedding $X\mapsto (X, -\sigma X)$. We let
$$
F(M) = \{m\in M\mid \dim U(\g) m <\infty\}.
$$
This is a $\g\times \g$-submodule of $M$. A $\g\times \g$-module $M$ 
is called a Harish-Chandra module if $F(M) = M$. We let $\HC$ denote
the category of Harish-Chandra modules.

\subsection{Constructions}

Let $M$ and $N$ be $\g$-modules. Then $\Hom_\C(M, N)$ and $(M\otimes_\C N)^*$
are $\g\times \g$-modules. We let
\begin{align*}
\L(M, N) &= F (\Hom_\C(M, N)) \\
\D(M, N) &= F( (M\otimes_\C N)^*)
\end{align*}
If $M\in \O$, then $\D(M, N) = \L(N, D M)$.

\subsection{Principal series modules in $\O_\lambda$}

Let $\lambda\in C$ and $\mu$ a dominant regular weight such  
that $\mu - \lambda\in \Z R$. Then 
$$
M\mapsto \L(M(\mu), M)
$$
defines an equivalence of $\O_\lambda$ with a subcategory $\tilde{\HC}$ 
of $\HC$. This result is due to Bernstein-Gelfand, Joseph, Enright (see 
Chapter 6 in \cite{Jantzen2}). 
The principal series modules in $\tilde{\HC}$ are
$$
M(x, y) = \D( M(y\cdot \lambda), M(x^{-1}\cdot\mu) )
$$
where $x, y\in W$. Via the above equivalence these can be viewed
as $\g$-modules in $\O_\lambda$. To stress this we sometimes use
the notation $M_\lambda(x, y)$.

\subsection{Twisted Verma properties}

In the following example and propositions we show
that principal series modules satisfy the properties of
twisted Verma modules.

\begin{example}
We have the following chain of isomorphisms ($\lambda$ and
$\mu$ as above)
\begin{align*}
M(x, e) &= \D( M(\lambda), M(x^{-1}\cdot \mu)) = \L(M(x^{-1}\cdot\mu), 
D M(\lambda)) \\
&= \L(M(x^{-1}\cdot \mu), M(\lambda)) = \L(M(\mu), M(x\cdot\lambda))
\end{align*}
where the last equality follows from (\cite{Jantzen2}, 7.23). This shows 
$M_\lambda(x, e) = M(x\cdot \lambda)$ and
that Property \ref{itemverma})
holds for principal series modules.
\end{example}

\begin{proposition}
Suppose that $\lambda, \bar{\lambda}\in C$. If $\lambda\in C^\circ$, then 
$$
T_\lambda^{\bar{\lambda}} M_\lambda(x, y) \cong
M_{\bar{\lambda}}(x, y)
$$ 
for every $x, y\in W$.
\end{proposition}
\begin{proof}
This follows from the corresponding property
$$
T_\lambda^{\bar{\lambda}} M (y\cdot \lambda) \cong M(y\cdot \bar{\lambda})
$$
for Verma modules
and the fact that the translation functor $T_\lambda^{\bar{\lambda}}$ 
becomes left translation of Harish-Chandra modules under the
equivalence $M\mapsto \L(M(\mu), M)$ (see \cite{Jantzen2}, 6.33).
\end{proof}
This verifies Property \ref{itemwall}). The following
proposition shows that Property \ref{itemshortexact}) holds.

\begin{proposition}
Let $\lambda\in C^\circ$. If $y s > y$, then
there is an exact sequence
$$
0\rightarrow M(x, y) \rightarrow \theta_s M(x, y) \rightarrow M(x, ys)\rightarrow 0,
$$
in $\O_\lambda$, where the first homomorphism is the adjunction map. 
\end{proposition}
\begin{proof}
This is (\cite{Irving}, Theorem 2.1).
\end{proof}
The following proposition is Property \ref{itemkashiwara}) verbatim.

\begin{proposition}
Let $\lambda\in C^\circ$.
Suppose that $x < x s $ and $s y > y$. Then we have
an isomorphism $M(x, y)\cong M(x s, s y)$ in $\O_\lambda$.
\end{proposition}
\begin{proof}
This is (\cite{Irving}, Theorem 4.4).
\end{proof}
By verifying the five properties of \S \ref{twistedvermaproperties} we have
proved that the principal series modules form a set of twisted
Verma modules by  Theorem \ref{theoremunique}.

\section{Local cohomology}
 
Let $G$ be a complex semisimple algebraic group with Lie algebra
$\g$, $T\subseteq B\subseteq G$ a maximal torus and a Borel
subgroup with Lie algebras $\h$ and $\b$ respectively. 
Let $X = G/B^{-}$
be the flag manifold of $G$, where $B^-$ is the Borel subgroup
opposite to $B$ and let $C(w)$ denote the $B$-orbit $Bw B^-/B^-$ in $X$.
Notice that $\codim\, C(w) = \ell(w)$.
A representation $M$ of $B^-$ induces a $G$-equivariant vector 
bundle $\L(M)$ on $X$. We let $X(B^-)=X(T)$ denote
the $1$-dimensional representations of $B^-$. Notice that
$X(T)$ can be identified with the integral weights in $\h^*$.
In general a $G$-linearized sheaf $\F$ of 
$\O_X$-modules is naturally a sheaf of $\hat{G}$-modules (where
$\hat{G}$ is the formal group of $G$) (\cite{Kempf}, Lemma 11.1) 
or equivalently a sheaf of $\Dist(G)\cong U(\g)$-modules. The local
cohomology group $\H^i_C(X, \F)$ has a natural $U(\g)$-module
structure for any locally closed subset $C\subseteq X$, where $i\geq 0$
(\cite{Kempf}, Lemma 11.1).

For a $B^-$-representation $M$ and a
locally closed subset $C\subseteq X$, we let $\H^i_C(M)$ denote
the $i$-th local cohomology group of $\L(M)$ with support
in $C$ with its natural $\g$-action.
By (\cite{Kempf}, Lemma 12.8) the
local cohomology groups are $\h$-diagonalizable and
$$
\ch \H^{\ell(w)}_{C(w)}(\lambda) = \ch M(w\cdot \lambda).
$$
This implies by Example \ref{exampleinO} that 
$\H^{\ell(w)}_{C(w)}(\lambda)\in \O$ (it belongs
in fact to the block $\O_\lambda$).

\subsection{Basic properties of local cohomology}

Local cohomology exists only in one degree in the following
sense.

\begin{proposition}
\label{proponedegree}
Let $V$ be a vector bundle on $X$ and $C$ an irreducible
affinely embedded locally closed subset of $X$ of codimension $\ell$. Then
$$
\H^i_C(X, V) = 0 \text{\ \ if\ } i\neq \ell.
$$
\end{proposition}
\begin{proof}
On the level of sheaves $\mathcal{H}^i_C(V)=0$ if $i\neq\ell$, since
$X$ is Cohen Macaulay and $C$ irreducible of codimension $\ell$. Now
one uses the local to global spectral sequence
$$
\H^p(X, \mathcal{H}^q_C(V))\implies \H^{p+q}_C(X, V)
$$
and the higher cohomology vanishing $\H^p(X,  \mathcal{H}^q_C(V)) = 0, 
p>0$, which follows from the assumption that $C$ is affinely
embedded, to deduce the result.
\end{proof}

\begin{proposition}
\label{propindenfor}
Let $V$ be a $B^-$ representation, $C$ a locally closed
subset of $X$ and $E$ a finite
dimensional $\g$-representation. Then there is an isomorphism 
$$
\H^i_C(V\otimes E) \cong \H^i_C(V)\otimes_\C E
$$
of $g$-modules for $i\geq 0$.
\end{proposition}
\begin{proof}
We may lift $E$ to a $G$-representation.
On the level of $G$-sheaves we have an isomorphism
$\L(V\otimes_\C E)\cong \L(V)\otimes_\C E$. This extends to
an isomorphism of $\hat{G}$-sheaves giving the desired result.
\end{proof}

\subsection{Principal series modules and local cohomology}

We emphasize the following important lemma.

\begin{lemma}[Kashiwara]
Let $\alpha\in S$ be a simple root, $w\in W$ and suppose that 
$\mu\in X(T)$ with $\<\mu, \alpha^\vee\> \geq -1$ and that
$w s_\alpha < w$. Then there is an isomorphism
$$
\H^{\ell(w)}_{C(w)}(\mu) \cong \H^{\ell(w)-1}_{C(w s_\alpha)}(
s_\alpha\cdot\mu)
$$
of $\g$-modules.
\end{lemma}
\begin{proof}
This is Lemma 3.6.6 in \cite{Kashiwara}.
\end{proof}
Fix $\lambda\in C^\circ\cap X(T)$. We will prove
that 
$$
\H^{\ell(x)}_{C(x)}(y\cdot \lambda)
$$
satisfies the properties of twisted Verma modules, thereby showing 
the isomorphism
$$
M(x, y)\cong \H^{\ell(x)}_{C(x)}(y\cdot \lambda)
$$
between principal series modules in $\O_\lambda$ and local cohomology.
Kashiwara's lemma is the key input for proving Property \ref{itemkashiwara}). 
In the 
above notation it states

\begin{lemma}[Kashiwara']
Let $\alpha\in S$ be a simple root and let $x, y\in W$, such that 
$x < x s_\alpha$ and $s_\alpha y > y$. Then there is an 
isomorphism
$$
\H^{\ell(x)}_{C(x)}(y\cdot\lambda) 
\cong \H^{\ell(x) +1}_{C(x s_\alpha)}(s_\alpha y\cdot\lambda)
$$
of $\g$-modules.
\end{lemma}
\noindent
The above lemma is the content of Property \ref{itemkashiwara}) for local
cohomology modules.

\begin{proposition}
There is an isomorphism
$$
\H^0_{C(e)}(\lambda) \cong D M(\lambda)
$$
of $\g$-modules for any (integral) weight $\lambda\in X(T)$.
\end{proposition}
\begin{proof}
This is Proposition 3.6.2 in \cite{Kashiwara}.
\end{proof}
The above proposition shows that Property \ref{itemverma}) holds for local
cohomology modules.

\subsection{Translation and local cohomology}

\begin{proposition}
\label{propfilt}
Let $0\rightarrow K\rightarrow V\rightarrow L\rightarrow 0$ be an exact
sequence of $B^-$ modules. Then we get an exact sequence
$$
0\rightarrow \H^i_{C(w)}(K)\rightarrow \H^i_{C(w)}(V) \rightarrow
\H^i_{C(w)}(L)\rightarrow 0
$$
of $\g$-modules for every $i\geq 0$ and $w\in W$.
\end{proposition}
\begin{proof}
This follows from the long exact sequence and 
Proposition \ref{proponedegree}.
\end{proof}
Let $\Pi(\eta)$ denote the weights in the finite dimensional
simple representation with extremal (integral) weight $\eta$.
We have the following special case of a well known lemma 
due to Jantzen (\cite{Jantzen1}, 2.9).

\begin{lemma}
\label{lemmajantzen}
Let $y\in W$ and $\lambda, \mu\in C$, where $\lambda\in C^\circ$ and $\mu-\lambda$ is integral. Then 
$$
W\cdot\mu \cap (y\cdot\lambda + \Pi(\mu-\lambda)) = \{y\cdot\mu\}.
$$
If $\Stab_W(\mu) = \{1, s\}$, then 
$$
W\cdot\lambda\cap(y\cdot\mu + \Pi(\lambda - \mu)) = \{y\cdot\lambda, 
y s\cdot\lambda\}.
$$
\end{lemma}

The following proposition shows that local cohomology modules
satisfy Properties 
\ref{itemshortexact}) and \ref{itemwall}) of twisted Verma modules with
the exception that one only has the short exact sequence in Property
\ref{itemshortexact}) (not knowing that the injection is the 
adjunction morphism). This unpleasant feature is resolved through
Lemma \ref{lemmaunique} and the construction of twisted Verma modules
using twisting functors (see Sections \ref{SectionTwistedVerma} and
\ref{SectionDerivedTwisting}). 

\begin{proposition}
Suppose that $y\in W$ and $\lambda, \mu\in C$, where $\lambda\in C^\circ$. 
Then
$$
T_\lambda^\mu \H^i_{C(w)}(y\cdot \lambda) = \H^i_{C(w)}(y\cdot \mu).
$$
If $\Stab_W(\mu)=\{1, s\}$ and $y s\cdot \lambda > y\cdot \lambda$, then
we have a short exact sequence
$$
0\rightarrow \H^i_{C(w)}(y\cdot\lambda)\rightarrow 
T_\mu^\lambda\H^i_{C(w)}(y\cdot\mu) \rightarrow \H^i_{C(w)}(y s\cdot\lambda)\rightarrow 0
$$
for every $i\geq 0$ and $w\in W$.
\end{proposition}
\begin{proof}
We use Proposition \ref{propindenfor}:
\begin{align*}
T_\lambda^\mu H^{\ell(w)}_{C(w)}(y\cdot\lambda) &= 
\pr_\mu (H^{\ell(w)}_{C(w)}(y\cdot\lambda)\otimes_\C E)\\
&= \pr_\mu \H^{\ell(w)}_{C(w)}(y\cdot \lambda\otimes E)
\end{align*}
where $E$ is the finite dimensional simple module with extremal
weight $\mu-\lambda$.
Observe that $\H^i_{C(w)}(\eta)\in \O_\eta$ for arbitrary $i\geq 0, 
w\in W$ and $\eta\in X(T)$.
Now take a $B^-$-filtration $N=N_0 \supseteq N_1 \supseteq \dots$ of 
$N=y\cdot\lambda\otimes E$, 
such that $N_i/N_{i+1} = \mu_i$ and $i < j \implies \mu_i \not< \mu_j$.
Then use Proposition \ref{propfilt} and Lemma \ref{lemmajantzen} to
get the desired result.
\end{proof}

\begin{remark}
Notice that we have proved the duality statement
$$
D\H^{\ell(w)}_{C(w)}(X, \L(\lambda))\cong
\H^{\ell(w w_0)}_{C(w w_0)}(X, \L(w_0 \cdot \lambda))
$$ 
of $\g$-modules
for arbitrary integral weights $\lambda$ and Schubert cells $C(w)$. This
follows from Corollary \ref{corollaryD}. 
\end{remark}

\section{Reformulation of formal properties of twisted Verma modules}

\label{SectionReform}

In this section we reformulate the properties in Section 2 describing a
family of twisted Verma modules. This is partly because we want to introduce
a new notation which is more natural in the setup in the following sections 
and partly because we want to generalize to the case of non-integral
weights. Of course, the principal series modules considered in Section 3 
also exist for non-integral weights (and in fact our definitions and results
in Section 3 immediately generalize to this case, see \cite{Irving}).

We fix an arbitrary weight $\lambda_0 \in \h^*$ and set $\Lambda = \lambda_0
+ \Z R \in \h^*/\Z R$. Then $R(\lambda_0) = \{ \alpha \in R \mid \;\langle \lambda_0, \alpha^\vee
\rangle \in \Z \}$ is a root system with corresponding Weyl group $W(\lambda_0) =
\{ w \in W \mid \; w(\lambda_0) - \lambda \in \Z R \}$. A weight $\lambda \in \Lambda$
is called dominant (respectively antidominant) if 
$\langle \lambda + \rho, \alpha^\vee \rangle \geq 0$ (respectively $\leq 0$)
for all $\alpha \in  R(\lambda_0) \cap R^+$. 

We define $\O_\Lambda$ to be the subcategory of $\O$ consisting of those
$M$ whose weights all belong to $\Lambda$.

\begin{definition}
\label{twistdef}
A family of twisted Verma modules in $\O_\Lambda$ is a collection of modules
$(M^w(\lambda))$ parametrized
by $\lambda\in\Lambda$ and $w\in W(\lambda_0)$ such that $M^w(\lambda) \in \O_\lambda$.
It is required to have the following properties
\begin{enumerate}[i)]
\item $M^e(\lambda') = M(\lambda')$ for some regular antidominant weight $\lambda'\in 
\Lambda$.
\item
Let $w, y, s\in W(\lambda_0)$, where $s$ is a simple reflection. If $ws > w$ and
$w^{-1} y < s w^{-1} y$ then we have an isomorphism
$M^w(y\cdot\lambda')\cong M^{w s}(y\cdot \lambda')$.

\item Let $w, y, s\in W(\lambda_0)$, where $s$ is a simple reflection. If 
$w^{-1} y > w^{-1} y s$ then we have a short exact sequence
$$
0\rightarrow M^w(y\cdot\lambda') \rightarrow \theta_s M^w(y\cdot\lambda')
\rightarrow M^w(ys\cdot\lambda')\rightarrow 0.
$$
\item For every antidominant weight $\mu\in \Lambda$ we have
$T_{\lambda'}^\mu M^w(\lambda) = M^w(\mu)$ for all $w\in W(\lambda_0), \lambda \in 
W(\lambda_0) \cdot \lambda'$.
\end{enumerate}
\end{definition}

In the next section we construct a family of twisted Verma modules and prove that
all its modules have $1$-dimensional endomorphism rings. Just as in 
Lemma 2.2 this shows that for any family of twisted Verma modules the first
homomorphism in the exact sequence appearing in Property iii) is (up to a nonzero
scalar) the adjunction morphism.
As in Section 2 this leads to the following results

\begin{theorem}
There is a unique family of twisted Verma modules in $\O_\Lambda$.
\end{theorem}

\begin{corollary}
\label{CorollaryDuality}
If $(M^w(\lambda))_{\lambda \in \Lambda, w \in W(\lambda_0)}$ is a family of
twisted Verma modules then $DM^w(\lambda) = M^{ww_0}(\lambda)$ for all $\lambda
\in \Lambda, w \in W(\lambda_0)$.
\end{corollary}

\begin{remark}
The correspondence between the above concept of a family of twisted Verma
modules and the previously considered one is given by
$$ M(x,y) = DM^x(xy\cdot \lambda').$$
It is straightforward to get the properties of $M(x,y)$ in 2.1 from the 
corresponding properties above.
\end{remark}
\section{Twisting functors}
In this section we consider the twisting functors introduced by Arkhipov
\cite{Arkhipov1}.

\subsection{The semiregular modules}
\label{Sectionsemireg}
Let $\g = \n^- \oplus \h \oplus \n^+$ be the triangular decomposition
of our semisimple complex Lie algebra $\g$ as in Section 1. Recall that $\b =
\h\oplus \n ^+$ is the
Borel subalgebra corresponding to $R^+$. We shall write 
$U = U(\g)$,
$N = \U(\n^-)$ and $B = \U(\b)$. 

The natural $\Z R$-grading on $\g$ (where elements in $\h$ 
have degree $0$ and elements in $\g_\alpha$ have degree $\alpha$,
$\alpha\in R$) gives rise to a grading on $U$, 
$$
U\cong \bigoplus_{\lambda\in\Z R} U_\lambda.
$$ 
Let $\hgt: \Z R\rightarrow \Z$ be the $\Z$-linear height
function with $\hgt(\alpha)=1$ for all simple roots $\alpha$. Then we
get a $\Z$-grading $U\cong \oplus_{n\in \Z} U_n$, where
$$
U_n = \bigoplus_{\hgt(\lambda)=n} U_\lambda,\,\,\,n\in \Z.
$$
Note that the subalgebra $N \subseteq U$ is negatively
graded with $N_0^-=\C$.

For $w \in W$ we consider the subalgebra $\n_w = \n^- \cap w^{-1}(\n^+)$ of
$\n^-$. The corresponding enveloping algebra $N_w = U(n_w)$ is then a
(negatively) graded subalgebra of $U$ with $(N_w)_0 = \C$. Note that $N_e = \C$
and $N_{w_0} = N$.

The (graded) dual of $N_w$ is $N^{*}_w = \oplus_{n\geq 0}
\Hom_\C((N_w)_n, \C)$. This is a $\Z$-graded bimodule over
$N_w $ with $(N^*_w)_n = \Hom_\C((N_w)_{-n} , \C), n\in \Z$. The left action 
of $N_w$ on $N_w^*$ is given by $x f : n\mapsto f(n x),
f\in N^{*}_w, x,n\in N_w$. The right
action  is defined similarly. 

Then we define the corresponding
semiregular module $S_w$ by 
$$
S_w = U\otimes_{N_w} N^{*}_w.$$

Clearly, $S_w$ is a left $U$-module and a right $N_w$-module. It is a non-trivial 
fact (see the theorem below) that $S_w$ is in fact a $U$-bimodule. To
state the precise result which gives this we first need a little more notation.

Let $e \in \n^-\setminus \{0\}$. Then we set $U_{(e)} = U \otimes_{\C[e]} \C[e, e^{-1}]$. In particular,
we shall consider the case where $e $ is equal to the Chevalley generator $e_{-\alpha}
\in \g_{-\alpha}$, $\alpha$ a simple root. Using this notation we can state

\begin{theorem}(Arkhipov \cite{Arkhipov1})
\begin{itemize}
\item [i)] 
\label{Ark1}
For each $ 0 \neq e \in \n^-$ we have that $U_{(e)}$ is an associative algebra
which contains $U$ as a subalgebra. We set $S_e = U_{(e)}/U$.

\item [ii)]
\label{Ark2} For each simple root $\alpha$ with corresponding simple reflection $s \in W$ 
we have an isomorphism of left $U$-modules $S_s \simeq S_{e_\alpha}$.

\item [iii)]
\label{Ark3} Let $w\in W$ and choose a filtration $\n_w = F^0 \supset F^1 \supset \cdots
\supset F^r \supset 0$ consisting of ideals $F^p \subset \n^-$ of codimension
$p , \; p = 0, 1, \cdots , r= l(w)$. If $e_p \in F^{p-1} \setminus F^p$ then
we have an isomorphism of $U$-bimodules
$$ S_w \simeq S_{e_1} \otimes _U  \cdots \otimes _U S_{e_r}.$$

\item [iv)]
\label{Ark4} (cf. Theorem 1.3 in \cite{WS}) For each $w \in W$ we have an isomorphism of right $U$-modules
$S_w \simeq N_w^*\otimes_{N_w}U$.
\end{itemize}
\end{theorem}

\subsection{The twisting functors on $\O$}
Let $\phi_w \in \Aut(\g)$ denote an automorphism corresponding to $w \in W$.
If $M$ is a $\g$-module we can conjugate the action of $\g$ on $M$ by $\phi_w$.
The module obtained in this way we shall denote $\phi_w(M)$. Note that if 
$\lambda \in \h^*$ then we have $ \phi_w(M)_{\lambda} = M_{w(\lambda)}$.

Following Arkhipov \cite{Arkhipov1} we define now a twisting functor $T_w$ on the
category of $\g$-modules by
$$ T_w M = \phi_w(S_w\otimes_UM).$$

\begin{remark} \label{remtwist}
\begin{itemize}
\item [i)] It is clear from the definition that $T_w$ is a right exact functor
for all $w \in W$.
\item [ii)] Theorem {\ref {Ark3}} iii) shows that we have $T_{ws} = T_w \circ T_s$
whenever $s$ is a simple reflection for which $ws > w$.

\end{itemize}
\end{remark}

We shall now consider the composite of the twisting functor with induction
from the subalgebra $B$.

Let $E$ be a left $B$-module  and set $T_w^{B} E = T_w(U\otimes_{B}E)$. 
Using Theorem \ref{Ark4} iv) and the fact that $U = N\otimes B$ we se that we
may identify $T_w^{B} E$ with $\phi_w(N_w^*\otimes_{N_w} N \otimes E)$ 
(as vector spaces and as $\h$-modules).
Here and elsewhere $\otimes$ without a subscript denotes tensor 
product over $\C$.

\begin{proposition} \label{propexact}
Let $w \in W$.
\begin{itemize}
\item[i)] The functor $T_w^{B}$ is exact.

\item[ii)]$ \ch T_wM(\lambda) = \ch T_w^{B} \lambda = \ch M(w\cdot \lambda)$ for all $\lambda \in \h^*$.
\end{itemize}
\end{proposition}
\begin{proof}  i) follows from the above by observing that $N$ is free over $N_w$. Also
ii) follows from the above identification via Theorem \ref{Ark1} i) and an easy
induction on $w$.
\end{proof}

As an immediate consequence of Proposition \ref{propexact} ii) we get
\begin{corollary}
\label{twistcor}
The functor $T_w$ restricts to a functor from $\O$ to $\O$ and it preserves $\O_\lambda$
for all $\lambda \in \h^*$.
\end{corollary}
\subsection{The tensor identity}
The following tensor identity will be important in the following
\begin{proposition}
\label{TIprop}
Let $w \in W$ and suppose $M$ and $V$ are $U$-modules with $V$ finite 
dimensional. Then we have a natural isomorphism $T_w(M\otimes V) \simeq
(T_w M)\otimes V$.

Likewise, if $E$ is a $B$-module then $T_w^{B} (E\otimes V) \simeq 
(T_w^{B} E) \otimes V$.
\end{proposition}
\begin{proof}
Note that by Remark \ref{remtwist} ii) we may reduce to the case where $w= s$ for some
simple reflection $s$. Since $\phi_s(V) \simeq V$ what we need to
prove is $S_s \otimes_U(M\otimes V) \simeq (S_s \otimes_U M)\otimes V$. But this
is clear from Theorem \ref{Ark1} i).

The last statement follows from the first by noting that the tensor identity for
induction ensures that we have an isomorphism $U\otimes_{B}(E\otimes V) \simeq
(U\otimes_{B}E)\otimes V$.
\end{proof}

\begin{corollary} \label{twistcom}
Let $w \in W$ and let $\lambda, \mu \in C$. Then
$T_w$ commutes with $T_\mu^\lambda$.

\end{corollary}

\begin{proof}
Let $M \in \O$ and write 
$M  = \bigoplus_\lambda \pr_\lambda M $. Since $T_w$
preserves $\O_\lambda$ (Corollary \ref{twistcor}) it follows that 
$T_w(\pr_\lambda M) = \pr_\lambda (T_w(M))$. Combining this with 
Proposition \ref{TIprop}
we get the statement.
\end{proof}

\subsection{Twisted Verma modules}
\label{SectionTwistedVerma}
Let $\lambda \in \h^*$ and $w \in W$. Then we define the twisted Verma module
$M^w(\lambda)$ by
$$M^w(\lambda) = T_wM(w^{-1} \cdot \lambda) = T_w^{B}(w^{-1} \cdot \lambda).$$

\begin{theorem} \label{thmtwist}
Let $\Lambda = \lambda_0 + \Z R \in \h^*/\Z R$. Then 
$(M^w(\lambda))_{\lambda \in \Lambda, w \in W(\lambda_0)}$
is a family of twisted Verma modules (in the sense of Definition \ref{twistdef}).
\end{theorem}
\begin{proof} 
Note that $M^w(\lambda) \in \O_\lambda$ by Corollary {\ref {twistcor}}.
We shall verify properties i), iii) and iv) in Definition 5.1 leaving ii) for later.

The above definition gives $M^e(\lambda) = M(\lambda)$ for all $\lambda$. So
property i) is certainly satisfied. Via 
Corollary 6.2 we see that property iii) is a consequence of the corresponding fact for ordinary
Verma modules.

Let now the notation and assumptions be as in Definition 5.1 iv). The well known effect of the
wall crossing functor $\theta_s$ on ordinary Verma modules gives the short exact sequence
$$0 \rightarrow M(w^{-1}y\cdot \lambda') \rightarrow \theta_sM(w^{-1}y\cdot \lambda')
\rightarrow M(w^{-1}ys\cdot \lambda' ) \rightarrow 0.$$
Applying the twisting functor $T_w$ to this sequence we get (using Corollary \ref{twistcom}
on the middle term)
$$0 \rightarrow M^w(y\cdot \lambda') \rightarrow \theta_sM^w(y\cdot \lambda')
\rightarrow M^w(ys\cdot \lambda' ) \rightarrow 0.$$
The exactness of this sequence comes from Proposition \ref{propexact} i).
\end{proof}

\subsection{The $\spl_2$-case}

Let $\g = \spl _2(\C)$ with the usual basis $\{f,h,e\}$. Then $\n^- = \C f, \h = \C h$ and 
$ \b = \h + \C e$. For $\lambda \in \h^* = \C$ the Verma module $M(\lambda)$ is simple
unless $\lambda \in \N$. On the other hand, when $\lambda \in \N$ we have an
exact sequence 
$$ 0 \rightarrow M(-\lambda - 2) \rightarrow M(\lambda) \rightarrow L(\lambda)
\rightarrow 0.$$
In this case there is only one non-trivial element $s$ in $W$. It is easy to see 
that $M^s(\lambda) = DM(\lambda)$. Hence $M^s(\lambda)$ is simple for $\lambda \notin
\N$ and for $\lambda \in \N$ we have an exact sequence
$$ 0 \rightarrow L(\lambda) \rightarrow M^s(\lambda) \rightarrow M^s(-\lambda -2)
\rightarrow 0.$$
Combining the above two sequences we get the following four term exact sequence (still assuming $\lambda \in \N$) 
$$ 0 \rightarrow M(-\lambda - 2) \rightarrow M(\lambda) \rightarrow M^s(\lambda)
\rightarrow M^s(-\lambda -2)\rightarrow 0.$$
Note that $M^s(-\lambda -2) = L(-\lambda -2) = M(-\lambda - 2)$.

\subsection{Twist and induction}
Returning to the general case we pick $\lambda \in \h^*$ and fix a simple root
$\alpha$. We denote by $\p_\alpha$ the minimal parabolic subalgebra of $\g$ containing 
$\b$ corresponding to $\alpha$. Then the $\p_\alpha$-Verma module with 
highest weight $\lambda$ is 
$$ M_\alpha(\lambda) = U(\p_\alpha)\otimes_{B} \lambda.$$

In analogy with the above we can also define an $s_\alpha$-twisted Verma module for $\p_\alpha$ 
with highest weight $\lambda$, namely
$$M_\alpha^{s_\alpha}(\lambda) = \phi_{s_\alpha}(U(\p_\alpha)\otimes_{N_{s_\alpha}}N_{s_\alpha}^*
\otimes_{U(\p_\alpha)}M_\alpha(s_\alpha \cdot \lambda))$$
Note that $N_{s_\alpha} = U(\g_{-\alpha}) = \C[e_{-\alpha}]$.

Considered as a module for the Levi subalgebra in $\p_\alpha$ we have that
$M_\alpha^{s_\alpha}(\lambda)$ is dual to $M_\alpha(\lambda)$. In particular we
get therefore as in the $\spl_2$-case:
\begin{lemma} \label{lemparabolic}
\begin{itemize}
\item[i)] If $\langle \lambda + \rho, \alpha^{\vee}\rangle \notin \N $ then we
have an isomorphism of $\p_\alpha$-modules $M_\alpha(\lambda) \simeq 
M_\alpha^{s_\alpha}(\lambda)$.

\item[ii)] If $\langle \lambda + \rho, \alpha^{\vee}\rangle \in \N $ then we have
an exact sequence of $\p_\alpha$-modules
$$0\rightarrow M_\alpha(s_\alpha \cdot \lambda) \rightarrow M_\alpha(\lambda)
\rightarrow M_\alpha^{s_\alpha}(\lambda) \rightarrow M^{s_\alpha}_\alpha(s_\alpha 
\cdot \lambda) \rightarrow 0.$$
The two extreme terms in this sequence are isomorphic.

\end{itemize} \end{lemma}

In the above notation we clearly have
$$M(\lambda) = U \otimes_{U(\p_\alpha)} M_\alpha(\lambda)$$
for all $\lambda \in \h^*$. We claim that we have a similar transitivity 
result for twisted Verma modules, namely
$$M^{s_\alpha}(\lambda) = U \otimes_{U(\p_\alpha)} M_\alpha^{s_\alpha}(\lambda).$$
This follows directly from the definitions: $M^{s_\alpha}(\lambda) =
\phi_{s_\alpha}(S_{s_\alpha}\otimes_UM(s_\alpha \cdot \lambda)) =
\phi_{s_\alpha}(U\otimes_{N_{s_\alpha}}N_{s_\alpha}^*\otimes_U U \otimes_{B} 
s_\alpha \cdot \lambda)) =
\phi_{s_\alpha}(U\otimes_{U(\p_\alpha)}U(\p_\alpha)\otimes_{N_{s_\alpha}}N_{s_\alpha}^*
\otimes_{B}s_\alpha \cdot \lambda) =
U\otimes_{U(\p_\alpha)}\phi_{s_\alpha}
(U(\p_\alpha)\otimes_{N_{s_\alpha}}N_{s_\alpha}^*
\otimes_{B}s_\alpha \cdot \lambda)
=U\otimes_{U(\p_\alpha)}M_\alpha^{s_\alpha}(\lambda).$

When we combine this with the results Lemma \ref{lemparabolic} we get
\begin{lemma} \label{lem2}
\begin{enumerate}
\item[i)] If $\langle \lambda + \rho, \alpha^\vee \rangle
\notin \N$ then we have an isomorphism of $\g$-modules 
$M(\lambda) \simeq M^{s_\alpha}(\lambda)$.

\item[ii)] If  $\langle \lambda + \rho, \alpha^\vee \rangle
\in \N$ then we have an exact sequence of $\g$-modules
$$0 \rightarrow M(s_\alpha \cdot \lambda) \rightarrow M(\lambda)
\rightarrow M^{s_\alpha}(\lambda) \rightarrow 
M^{s_\alpha}(s_\alpha \cdot \lambda) \rightarrow 0.$$
\end{enumerate}
\end{lemma}

Finally, let $w \in W$ and suppose $ws_\alpha > w$. 
Then $\beta = w(\alpha) \in \R^+$. By Remark \ref{remtwist} ii) we have
$T_{ws_\alpha} = T_w \circ T_{s_\alpha}$ and hence by applying $T_w$
the results in the above lemma we find (replacing $\lambda$ by $w^{-1}\cdot \lambda)$

\begin{proposition}\begin{enumerate}
\item[i)] If $\langle \lambda + \rho, \beta^\vee \rangle \notin \N$
then $M^w(\lambda) \simeq M^{ws_\alpha}(\lambda)$.

\item[ii)] If $\langle \lambda + \rho, \beta^\vee \rangle \in \N$ then
we have an exact sequence of $\g$-modules 
$$0 \rightarrow M^w(s_\beta \cdot \lambda) \rightarrow M^w(\lambda)
\rightarrow M^{ws_\alpha}(\lambda) \rightarrow 
M^{ws_\alpha}(s_\beta \cdot \lambda) \rightarrow 0.$$
\end{enumerate}
\end{proposition}
The fact that $T_w$ preserves the exactness of the sequence 
in Lemma \ref{lem2} is a consequence of Proposition \ref{propexact} i).

\begin{remark}
\begin{itemize}
\item[i)] Let the notation and assumptions be as in Definition 5.1 ii)
and let $\alpha$ be the simple root with reflection $s$. Set $\beta =
w(\alpha)$. Then the assumption $w^{-1}y < sw^{-1}y $ is equivalent to
$w^{-1}y \cdot \lambda' < sw^{-1}y \cdot \lambda' $, i.e 
$\langle y(\lambda' + \rho), \beta^\vee\rangle = 
\langle w^{-1}y(\lambda' + \rho), \alpha^\vee \rangle < 0$.
Hence i) in this proposition gives $ M^w(y \cdot \lambda') \simeq M^{ws_\alpha}(y \cdot \lambda')$.
We have thus completed the proof of Theorem \ref{thmtwist}.

\item[ii)] It is possible to derive the four term exact sequence in this proposition
from the formal properties of a family of twisted Verma modules, see \cite{Irving}.
Here we have taken the opportunity to derive it directly from the $\spl _2$-case.
The sequence is sometimes called the Zelebenko-Duflo-Joseph sequence.
\end{itemize}
\end{remark}
\subsection{Endomorphisms}
\label{SectionDerivedTwisting}
The twisting functor $T_w$ on $\O$ extends to the functor $LT_w$ on the 
derived category $\D^b(\O)$. A key property of this functor is

\begin{proposition} (Arkhipov \cite{Arkhipov2})
Let $w \in W$. The derived functor $LT_w$ is an autoequivalence on
$D^b(\O)$.
\end{proposition}

\begin{corollary}
\label{CorollaryConstantEndomorphisms}
Let $\lambda \in \Lambda$ and $w \in W(\lambda_0)$. Then $\End_\O(M^w(\lambda)) 
= \C$.
\end{corollary}

\begin{proof}
The result is well known for ordinary Verma modules. It then follows for twisted
Verma modules by the above proposition. In fact, we get
$\End_\O(M^w(\lambda)) = \End_O(T_wM(w^{-1} \cdot \lambda)) = 
\End_{\D^b(\O)}(LT_wM(w^{-1}\cdot \lambda)) = \End_{\D^b(\O)}(M(w^{-1}\cdot \lambda))
=\End_\O(M(w^{-1}\cdot \lambda)) = \C$.
\end{proof}

\begin{remark}
This corollary is essential for our arguments in Section 2 proving the uniqueness of 
a family of twisted Verma modules, see Lemma \ref{lemmaunique} . More precisely, we need
for $\lambda$ regular and $s$ a simple reflection that
 $\Hom_\g(M^w(\lambda), \theta_sM^w(\lambda)) = \C$. But this Hom-space equals
$$\End_\g(T_\lambda^\mu M^w(\lambda)) = \End_\g (M^w(\mu))$$ by property iv) in
Definition \ref{twistdef}. In other words we need the corollary for semi-regular
weights. (In the Harish-Chandra module situation this is exactly the case 
handled by C. Stroppel, see Remark \ref{remStroppel}). 
\end{remark}

\section{Filtrations and sum formulae}

In this section we shall show that the twisted Verma modules considered in
the previous sections have Jantzen type filtrations and we shall give the
corresponding sum formulae.

\subsection{Deformations}
Let $A = \C[X]_{(X)}$ be the localization of the polynomial ring $\C[X]$ in
the maximal ideal generated by $X$. We shall then consider the Lie algebra $\g_A = 
\g\otimes_\C A$ over $A$. Similarly, we set $\h_A = \h\otimes_\C A,
\b_A = \b\otimes_\C A,$ etc.

If $\lambda \in \h^*_A = \Hom_A(\h_A, A) \simeq \h^*\otimes_\C A$ we have
a Verma module over $A$ \; $M_A(\lambda) = U_A\otimes_{B_A} \lambda$ where $U_A =
U(\g_A)$ and $B_A = U(\b_A)$. Note that $M_A(\lambda)$ is free over $A$
and that if $A\rightarrow \C$ is the specialization which takes $X$ into $0$
then we have $M_A(\lambda)\otimes_A\C \simeq M(\bar \lambda)$. Here 
$\bar \lambda = \lambda \otimes 1 \in \h^*_A\otimes_A\C = \h^*$.

The twisting functors $T_w, \; w\in W$ may also be defined over $A$. We just
extend scalars from $\C$ to $A$, i.e if $M$ is a $U_A$-module we set 
$T_w M = \phi_w(S_w^A \otimes _{U_A} M)$ with $S_w^A = S_w\otimes_\C A$.

In particular, this allows us to define twisted Verma modules over $A$
$$ M_A^w(\lambda) = T_w(M_A(w^{-1} \cdot \lambda)),$$
$\lambda \in \h_A^*, \; w \in W$. These modules specialize to the twisted Verma 
modules $M^w(\bar \lambda)$ considered in the previous section.

\subsection{The $\spl_2$-case revisited}
Consider again $\g = \spl_2(\C)$. Then the Verma module $M_A(\lambda)$ for
$\g_A = \spl_2(A)$ with highest weight $\lambda \in h^*_A = A$ has basis 
$\{v_0, v_1, \dots\}$ with
$v_0$  a generator for $M_A(\lambda)_\lambda$ and $v_i = f^{(i)} v_0, 
i\geq 0$. Here $f^{(i)} = f^{i}/i!\in U_A$. The action of
$\g_A$ on $M_A(\lambda)$ is given by
$$
h v_i = (\lambda - 2i) v_i,\,\,\, f v_i = (i+1) v_{i+1},\,\,\,
e v_i = (\lambda + 1 - i) v_{i-1}, i\geq 0.
\label{eq:10}
$$
for all $ i\geq 0$ (we set $v_{-1}=0$).

The twisted Verma module $M_A^s(\lambda)$ is
dual to $M_A(\lambda)$. So if $\{v_0^*, v_1^*, \dots\}$ denotes
the dual basis  it is immediate to check that the linear
map $\varphi_\lambda : M_A(\lambda) \rightarrow M^s_A(\lambda)$ given
by
$$
\varphi_\lambda(v_i) =  \binom{\lambda}{i} v_i^*, i= 0, 1, \dots
$$
is a $\g_A$-homomorphism which generates $\Hom_{\g_A}(M_A(\lambda), 
M_a^s(\lambda))\cong A$.

Note that if $\lambda \notin \N$ then the elements $ \binom{\lambda}{i} \in A$
are units in $A$ for all $i$. Hence in this case
$\varphi_\lambda$ is an isomorphism.
Define now
$\psi_\lambda: M^s(\lambda)
\rightarrow M(\lambda)$ to be the inverse of $\varphi_\lambda$ when $\lambda \notin \N$
and 
$$
\psi_\lambda(v_i^*)=
\begin{cases}
0, &\text{ if } i\leq \lambda\\
(-1)^i\binom{i}{i-\lambda-1}v_i, &\text{ if } i\geq \lambda+1
\end{cases}
$$
when $\lambda\in\N$. 

It is easy to check that $\psi_\lambda$ is a generator
of $\Hom_{\g_A}(M^s_A(\lambda), M_A(\lambda))$.
When we pass to the specialization $A \rightarrow \C$ we have of course still that
$\psi_{\bar \lambda}$ is an isomorphism when the specialization  $\bar \lambda$
of $\lambda$ is not in $\N$. If $\bar \lambda \in \N$
we get in analogy with the sequence involving $\varphi_{\bar \lambda}$ a four term exact sequence
$$0\rightarrow M(\bar \lambda)/M(-\bar\lambda - 2) \rightarrow M^s(\bar \lambda)
\rightarrow M(\bar \lambda) \rightarrow  M(\bar \lambda)/M(-\bar\lambda - 2)
\rightarrow 0$$
where the middle homomorphism is $\psi_{\bar \lambda}$.

Now we fix $\lambda \in \h^*$
and look at the character $\lambda + X \in \h^*_A$.
Then we have that
$$
M_A(\lambda+X)\stackrel{\varphi_\lambda}{\rightarrow} M_A^s(\lambda+X)
$$
is an isomorphism if $\lambda\not\in\N$ and fits into the exact sequence
$$
0\rightarrow M_A(\lambda + X) \stackrel{\varphi_\lambda}{\rightarrow}
M^s_A(\lambda+X)\rightarrow M(-\lambda-2)\rightarrow 0.
$$
if $\lambda\in\N$.

In fact, the first claim is a special case of the situation dealt with above. To 
verify the second statement we assume $\lambda\in \N$. Then we see that 
$\binom{\lambda + X}{i}$ is a unit in $A$ only when $0\leq i\leq\lambda$ 
whereas $\binom{\lambda + X}{i} = u_i X$ for some
unit $u_i\in A$ when $i>\lambda$. Moreover, we may
identify $M_A(-\lambda-2+X)/ X M_A(-\lambda-2+X)$ with
$M(-\lambda-2)$ and we get a surjection $M^s_A(\lambda+X)
\rightarrow M(-\lambda-2)$ by sending $v_i^*\mapsto (-1)^i v_{i-\lambda-1}$,
$i>\lambda$, $v_i^*\mapsto 0, i\leq\lambda$. It is now easy to check that 
this leads to an exact sequence as claimed.

Similarly, just as over $\C$ we have for all $\lambda$ a natural homomorphism
$\psi_\lambda: M^s_A(\lambda+X)\rightarrow M_A(\lambda+X)$. When
$\lambda\not\in\N$ this is the inverse of $\varphi_\lambda$ and for
$\lambda\in\N$ we have the exact sequence
$$
0\rightarrow M^s_A(\lambda+X)\stackrel {\psi_\lambda}{\rightarrow}
M(\lambda+X)\rightarrow M(\lambda)/M(-\lambda-2)\rightarrow 0.
$$

\subsection{The general case}

Consider now $\g$ general. For each simple root $\alpha$ the results
above transform easily into statements
about $\p_\alpha$-modules. So we may
proceed exactly as in Section 6 to obtain
the following results.

\begin{proposition} \label{prophom}
Let $\lambda\in\h^*$ and consider $\lambda+X\rho\in\h_A^*$.
Suppose $w\in W$ and $\alpha$ is a simple root with $w s_\alpha
> w$. Set $\beta = w(\alpha)$. Then $\Hom_{\g_A}(M^w_A(\lambda + X\rho),
M^{w s_\alpha}_A(\lambda + X\rho))\cong A\cong
\Hom_{\g_A}(M^{w s_\alpha}_A(\lambda + X\rho), M^{w}_A(\lambda + X\rho))$.
Moreover, if $\varphi_\lambda^w$ and $\psi_\lambda^w$ denote
generators of these $\Hom$-spaces. Then we have
\begin{itemize}
\item[i)]
If $\<\lambda + \rho, \beta^\vee\>\not\in \N$ then $\varphi_\lambda^w$ and
$\psi_\lambda^w$ are isomorphisms with $\varphi_\lambda^w =
(\psi_\lambda^w)^{-1}$.
\item[ii)]
If $\<\lambda+\rho, \beta^\vee\>\in\N$ then $\varphi_\lambda^w$ and
$\psi_\lambda^w$ fit into the exact sequences
$$
0\rightarrow M_A^w(\lambda+X\rho)\stackrel{\varphi_\lambda^w}{\rightarrow}
M_A^{w s_\alpha}(\lambda + X\rho)\rightarrow M^w(s_\beta\cdot\lambda)
\rightarrow 0
$$
and
$$
0\rightarrow M_A^{w s_\alpha}(\lambda+X\rho)\stackrel{\psi_\lambda^w}{\rightarrow}
M_A^{w}(\lambda + X\rho)\rightarrow M^w(\lambda)/M^w(s_\beta\cdot
\lambda)\rightarrow 0,
$$
respectively.
\end{itemize}
\end{proposition}

Fix now $\lambda\in\h^*$ and $w\in W$. Choose a reduced expression
for $w_0$, $w_0 = s_1 s_2 \dots s_N$ such that 
$w = s_n s_{n-1} \dots s_1$. Let $\alpha_{i_j}$ denote
the simple root corresponding to $s_j$. Set
$$
\beta_j = 
\begin{cases}
-w s_1 s_2 \dots s_{j-1}(\alpha_{i_j}), &\text{if } j\leq n\\
w s_1 s_2 \dots s_{j-1}(\alpha_{i_j}), &\text{if } j> n
\end{cases}
$$
Then $\{\beta_1, \beta_2, \dots, \beta_N\} = R^+$. If we set
$R^+(w) = \{\beta\in R^+\mid w^{-1}(\beta)\in R^-\}$, then
$\{\beta_1, \beta_2, \dots, \beta_n\} = R^+(w)$. We shall also write
$R^+(\lambda) = \{\beta\in R^+\mid \<\lambda + \rho, \beta^\vee\>
\in \N\}$.

Consider the composite $\Phi^w(\lambda)$
$$
M_A^w(\lambda + X\rho)\rightarrow
M_A^{w s_1}(\lambda+X\rho)\rightarrow
M_A^{w s_1 s_2}(\lambda+X\rho)\rightarrow
\dots\rightarrow
M_A^{w w_0}(\lambda+X\rho)
$$
where for each $j=1, \dots, N$ the homomorphism
$M_A^{w s_1 s_2 \dots s_{j-1}}(\lambda + X\rho)
\stackrel{\varphi_j^w(\lambda)}{\rightarrow} M_A^{w s_1 s_2 \dots s_j}
(\lambda + X\rho)$ is a generator of its $\Hom$-space (see
Proposition \ref{prophom}). Then we may define the
Jantzen filtration of $M_A^w(\lambda+X\rho)$ by
$$
M_A^w(\lambda+X\rho)^j=
\{m\in M_A^w(\lambda+X\rho)\mid
\Phi^w(\lambda)(m)\in X^j M_A^{w w_0}(\lambda+X\rho)\}.
$$
Taking the images in $M^w(\lambda) = 
M_A^w(\lambda + X\rho)/X M_A^w(\lambda+X\rho)$ we obtain
the Jantzen filtration $M^w(\lambda)^0\supseteq M^w(\lambda)^1\supseteq
\dots$ of $M^w(\lambda)$.

These filtrations also filter the weight spaces of
$M_A^w(\lambda+X\rho)$ and $M^w(\lambda)$. Note that
for any $\mu\in \h^*$ the weight space $M^w_A(\lambda)_{\mu+X\rho}$
is a finitely generated free $A$-module (of rank equal to
$\dim_\C M^w(\lambda)_\mu = \dim_\C M(\lambda)_\mu$). Standard
arguments (see e.g 5.1 in \cite{Jantzen1}) tell us that
$$
\sum_{j\geq 1} \dim M^w(\lambda)_\mu^j = \nu_X(\det(\Phi^w(\lambda)_\mu)).
$$
(Here and elsewhere the index $\mu$ on a homomorphism means the
restriction of the homomorphism to the $\mu+X\rho$ weight space and
$\nu_X$ is the $X$-adic valuation).

Clearly, the right hand side of this equation equals 
$$
\sum_{j=1}^n \nu_X(\det(\varphi_j^w(\lambda)_\mu)) = 
\sum_{j=1}^n \ell_X(\Coker(\varphi^w_j(\lambda)_\mu))
$$ 
where $\ell_X$ denotes length of a module. Observe that by
Proposition \ref{prophom} i) we have that $\varphi_j^w(\lambda)$
is an isomorphism when $\beta_j\not\in R^+(\lambda)$. By 
Proposition \ref{prophom} ii) we have for $\beta_j\in R^+(\lambda)$
$$
\ell_X(\Coker(\varphi_j^w(\lambda)_\mu)) =
\begin{cases}
\dim M(\lambda) - \dim M(s_{\beta_j}\cdot \lambda)_\mu,
&\text{if } j\leq n\\
\dim M(s_{\beta_j}\cdot \lambda)_\mu, &\text{if } j > n.
\end{cases}
$$
Hence we have proved
\begin{theorem} \label{thmsum}
Let $\lambda, w$ be as above. Then $M^w(\lambda)$ has a Jantzen
filtration
$$
M^w(\lambda) = M^w(\lambda)^0\supseteq M^w(\lambda)^1\supseteq\dots
$$
such that 
$
M^w(\lambda)/ M^w(\lambda)^1 \cong \Im \Phi^w(\lambda) \subseteq
M^{w w_0}(\lambda)
$
and
\begin{align*}
\sum_{j\geq 1} \ch M^w(\lambda)^j = 
&\sum_{\beta\in R^+(\lambda)\cap R^+(w)} (\ch M(\lambda) -
\ch M(s_\beta\cdot \lambda))  + \\
&\sum_{\beta\in R^+(\lambda)\setminus R^+(w)} \ch M(s_\beta\cdot \lambda).
\end{align*}
\end{theorem}

\begin{remark}
\begin{itemize}
\item[i)]
For $w=e$ we recover the usual Jantzen filtration and sum formula 
for the ordinary Verma module $M^e(\lambda) = M(\lambda)$ (Note
that in this case $R^+(w) = \emptyset$).
\item[ii)]
When reformulated using the notation from Sections 2--4 for twisted Verma modules
the theorem reads as follows:

Let $x, y \in W$. Then $M(x,y) = M_\lambda(x,y)$ (with $\lambda$ a regular 
integral and antidominant weight)
has a Jantzen filtration
$$M(x,y) = M(x,y)^0\supseteq M(x,y)^1\supseteq\dots
$$
such that
$
M(x,y)/ M(x,y)^1 \cong \Im (M(x,y) \rightarrow M(xw_0,w_0y))
$
and
\begin{align*}
\sum_{j\geq 1} \ch M(x,y)^j = 
&\sum_{\beta\in R^+(xy)\setminus R^+(x)} (\ch M(xy \cdot \lambda) -
\ch M(s_\beta xy\cdot \lambda)) + \\
&\sum_{\beta\in R^+(xy)\cap R^+(x)} \ch M(s_\beta xy\cdot \lambda).
\end{align*}
\end{itemize}
\end{remark}

\subsection{The $B_2$-case}
\begin{example}
Below we have listed all the twisted Verma modules with integral highest weights
when the root system is $B_2$ together with their Jantzen filtrations. Since
this is a multiplicity free case the sum formula in Theorem \ref{thmsum} 
completely determines the filtration. A simple module listed in the $i$-th 
row means that it occurs in the $i$-th layer of the filtration. In some cases 
a $0$ occurs in the $0$-th row. This means that the corresponding layer is $0$.

Choose an integral regular antidominant weight $\lambda$ and write $M^w(y)$ short
for $M^w(y \cdot \lambda)$. Also write $L(x) = L(x \cdot \lambda)$. Let $s$
(respectively $t$) be the simple reflections corresponding to the short 
(respectively long) simple root. Then $W = \{e, s, t, st, ts, sts, tst, w_0 \}$.

Recall that $M^w(y) = DM^{ww_0}(y)$ for all $w,y \in W$. Therefore we have only 
listed half the twisted Verma modules. The others (and their Jantzen filtrations)
are then obtained by dualizing. In the list below the twisted Verma
modules are itemized  according to their highest weight.

\begin{itemize}
\item[\fbox{$\lambda$}] $$M^w(e) = L(e)$$
\item[\fbox{$s\cdot\lambda$}] 
$$
M^{tst}(s) = M^{ts}(s) = M^t(s) = M^e(s) = 
\begin{matrix} L(s) \\ L(e) \end{matrix}
$$
\item[\fbox{$t\cdot\lambda$}]
$$
M^{tst}(t) = M^{st}(t) = M^s(t) = M^e(t) = 
\begin{matrix} L(t) \\ L(e) \end{matrix}
$$
\item[\fbox{$st\cdot\lambda$}]
\begin{align*}
&M^{ts}(st) = M^t(st) = M^e(st) = 
\begin{matrix}
L(st)\\
L(s)\,\,\,L(t) \\
L(e) 
\end{matrix}\\
&\\
&M^s(st) = 
\begin{matrix}
L(t)\\
L(e)\,\,\,L(st) \\
L(s)
\end{matrix}
\end{align*}
\item[\fbox{$ts\cdot\lambda$}]
\begin{align*}
&M^{st}(ts) = M^s(ts) = M^e(ts) = 
\begin{matrix}
L(ts)\\
L(s)\,\,\,L(t) \\
L(e) 
\end{matrix}\\
&\\
&M^t(ts) = 
\begin{matrix}
L(s)\\
L(e)\,\,\,L(ts) \\
L(t)
\end{matrix}
\end{align*}
\item[\fbox{$sts\cdot\lambda$}]
\begin{align*}
&M^{t}(sts) = M^e(sts) =  
\begin{matrix}
L(sts)\\
L(st)\,\,\,L(ts) \\
L(s)\,\,\,L(t) \\
L(e) 
\end{matrix}\\
&\\
&M^s(sts) = 
\begin{matrix}
L(ts)\\
L(s)\,\,\,L(t)\,\,\,L(sts) \\
L(e)\,\,\,L(st)
\end{matrix}\\
&\\
&M^{st}(sts) =
\begin{matrix}
0\\
L(e)\,\,\,L(ts)\\
L(s)\,\,\,L(t)\,\,\,L(sts)\\
L(st)
\end{matrix}
\end{align*}
\item[\fbox{$tst\cdot\lambda$}]
\begin{align*}
&M^{s}(tst) = M^e(tst) =  
\begin{matrix}
L(tst)\\
L(st)\,\,\,L(ts) \\
L(s)\,\,\,L(t) \\
L(e) 
\end{matrix}\\
&\\
&M^t(tst) = 
\begin{matrix}
L(st)\\
L(s)\,\,\,L(t)\,\,\,L(tst) \\
L(e)\,\,\,L(ts)
\end{matrix}\\
&\\
&M^{ts}(tst) =
\begin{matrix}
0\\
L(e)\,\,\,L(st)\\
L(s)\,\,\,L(t)\,\,\,L(tst)\\
L(ts)
\end{matrix}
\end{align*}
\item[\fbox{$w_0\cdot\lambda$}]
\begin{align*}
&M^{e}(w_0) =   
\begin{matrix}
L(w_0)\\
L(sts)\,\,\,L(tst) \\
L(st)\,\,\,L(ts) \\
L(s)\,\,\,L(t)\\
L(e) 
\end{matrix}\\
&\\
&M^s(w_0) = 
\begin{matrix}
L(tst)\\
L(st)\,\,\,L(ts)\,\,\,L(w_0) \\
L(s)\,\,\,L(t)\,\,\,L(sts)\\
L(e)
\end{matrix}\\
&\\
&M^{t}(w_0) =
\begin{matrix}
L(sts)\\
L(st)\,\,\,L(ts)\,\,\,L(w_0)\\
L(s)\,\,\,L(t)\,\,\,L(tst)\\
L(e)
\end{matrix}\\
&\\
&M^{st}(w_0) =
\begin{matrix}
0\\
L(t)\,\,\,L(tst)\\
L(e)\,\,\,L(st)\,\,\,L(ts)\,\,\,L(w_0)\\
L(s)\,\,\,L(sts)
\end{matrix}
\end{align*}
\end{itemize}
\end{example}
\begin{remark}
Recall that the Jantzen filtration of an ordinary Verma module is its unique
Loewy series, see \cite{Irving2}. In particular, the radical series of the Verma modules can
be read off from the above list and we have therefore a determination 
of all extensions between simple 
modules. Using this it is easy to see that there are twisted Verma modules
which do not have simple heads. For instance, both $L(e)$ and $L(ts)$ are 
quotients of $M^{st}(sts)$. Likewise, both $L(e)$ and $L(sts)$ are submodules
of $M^s(w_0)$ (this example of a non-rigid twisted Verma module was pointed
out to us by C. Stroppel). It is also seen that $M^{st}(w_0)$ has non-simple 
head and socle. 

The $0$ occurring in the $0$-th row for a module $M$ in the list means that the
composite $M \rightarrow DM$ (see 7.3) is zero. Nevertheless, the space
$\Hom_\g(M, DM)$ may be non-zero. For instance, one may check that 
$$\Hom_\g(M^{st}(w_0), DM^{st}(w_0))$$ is $2$-dimensional.
\end{remark}

\begin{remark}
Using that $LT_w$ is an autoequivalence of the bounded derived
category $D^b(\O)$ one may prove that (similar to the proof
of Corollary \ref{CorollaryConstantEndomorphisms})
$$
\Hom_\g(M^w(\lambda), M^{w s}(\lambda)) = \C
$$
where $w < w s$. Let $\varphi$ be a generator of this $\Hom$-space.
It seems reasonable to expect that $\varphi$ is well behaved 
with respect to the Jantzen filtration in the sense that 
$$
\varphi(M^w(\lambda)^j) \subseteq M^{w s}(\lambda)^{j+1}.
$$ 
One may prove that $\Soc M^{w_0 s_\alpha}(\lambda) = L(s_\alpha\cdot
\lambda)$, where $\alpha$ is a simple root and $\<\lambda+\rho, \alpha^\vee\>
\geq 0$.
If $\varphi$
respects the Jantzen filtration as above this leads to new
and perhaps simpler proofs of non-vanishing $\Ext^1$-groups
between certain neighboring simple modules.
\end{remark}

\end{document}